\documentclass[11pt]{amsart}
\usepackage{amsfonts,amssymb,amscd,amsmath,enumerate,verbatim,calc}

\newcommand{\m}{\mathfrak{m} }

\newcommand{\q}{\mathfrak{q} }

\def\boldy{\mbox{\boldmath $y$}}

\newcommand{\x}{ \mathbf{x}}
\newcommand{\y}{ \mathbf{y}}
\newcommand{\X}{ \mathbf{X}}

\newcommand{\ann}{\operatorname{ann}}

\theoremstyle{plain}
\newtheorem{thm}{Theorem}

\theoremstyle{definition}

\theoremstyle{remark}

\begin{document}

\title [partial Euler characteristics  ]{A short note on the non-negativity of partial Euler 
characteristics }
\author{Tony ~J.~Puthenpurakal}
\date{\today}
\address{Department of Mathematics, IIT Bombay, Powai, Mumbai 400 076}

\email{tputhen@math.iitb.ac.in}
 \begin{abstract}
Let 
$(A,\mathfrak{m})$ 
be a Noetherian local ring, $M$ 
a finite $A$-module and  $x_1,\ldots,x_n\in \m$ such that 
$\lambda (M/\x M)$ is finite.
Serre \cite[Appendix 2]{LoA} proved that all partial Euler 
characteristics  of $M$ with respect to $\x$ is non-negative.
This fact is easy to show when $A$ contains a field \cite[4.7.12]{BH}.
We give an elementary proof of Serre's result when $A$ does not contain a
field. 
\end{abstract}

 \maketitle

Let $(A,\mathfrak{m})$  be a Noetherian local ring and $M$ a finite 
$A$-module. Let $x_1,\ldots,x_n\in\m $ be a 
{\em multiplicity system} of $M$ i.e. $\lambda(M/(\x)M)$ is 
finite. (Here $\lambda (-)$ denotes length). Let $K(\x,M)$ 
be the Koszul complex of $\x$ with coefficients in $M$ and 
let $H_\bullet (\x,M)$ be its homology. Note that $H_\bullet 
(\x,M)$ has finite length.
One defines for all $j\geq 0$ the {\em partial Euler 
characteristics}
$$
\chi_j (\x,M)=\sum_{i\geq j} 
(-1)^{i-j}\lambda(H_i(\x,M))
$$
of $M$ with respect to $\x$. Serre showed all the partial 
Euler characteristics are non-negative. It is well known that 
$\chi_0(\x ,M)$ is either zero or the multiplicity of $M$ with 
respect to the ideal $(x_1,\ldots,x_n)$. It is also easy to see 
that $\chi_1(\x,M)$ is non-negative, \cite[4.7.10]{BH}.  The non-negativity of 
$\chi_j(\x,M)$ for $j\geq 2$ can be easily proved if $A$ 
contains a field,\cite[4.7.12]{BH}. In this short note we give an elementary proof 
of Serre's Theorem when $A$ does not contain a field. 
\begin{thm}
Let $(A,\m)$ be a Noetherian local ring, $A$ not containing a field. Let  $M$ be  a finite
 $A$-module 
and $x_1,\ldots, x_n\in \m$ a multiplicity system of $M$. Then
$$
\chi_j(\x,M)\geq 0\quad \mbox{for each } j\geq 0.
$$
\end{thm}
\begin{proof}
We may assume that $A$ is complete. To prove the theorem we 
construct a local Noetherian ring $(B,\m$) with a local 
homomorphism $\varphi :B\rightarrow A$, and $y_1,\ldots,y_n\in 
\m$ such that
\begin{enumerate}
\item
 $\varphi (y_i)=x_i$.
\item
 $M$ becomes a finite $B$-module (via $\varphi$).
\item
 $y_1,\ldots,y_n$ is a regular sequence and a s.o.p of 
$B$.
\end{enumerate}
Since $K(\y,M)\simeq K(\x,M)$ (as $B$-modules), we have 
$H_\bullet (\y,M) \cong H_\bullet (\x,M)  $
and so $\chi_j(\y,M)=\chi_j(\x,M)$ for each $j \geq 0$.

Suppose we have constructed $B$ as above.  The result then follows on similar lines
as  in \cite[4.7.12]{BH}. We give the proof here for the readers convenience.
 We prove the result by induction 
on $j$. For $j =0,1$ the result is already known. Let $j > 1$ and consider 
an exact sequence
$$
0\rightarrow U\rightarrow F\rightarrow M\rightarrow 0$$
where $F$ is a finite free $B$-module. Since $\boldy$ is $B$-regular we 
have $H_i(\y,F) = 0$ for $i >0$. Therefore
$ H_i(\y,M) \simeq H_{i-1}(\y,U) \quad\mbox{for all}\ i > 1.$
This yields $\chi_j (\y,M)=\chi_{j-1}(\y, U)$ and the proof is 
complete by induction hypothesis.

\emph{Construction of $B$}

Since $A$ is complete there exists a DVR, $(R,\rho)$ and a ring homomorphism
$\varphi : R\rightarrow A$ which induces an isomorphism $R/\rho R\rightarrow 
A/\m$.
Set $S= R[[ X_1,\ldots, X_n]]$ and let $\q$ be its maximal ideal and 
consider the natural ring map $\phi : S\rightarrow A$, with $ \phi (X_i)=x_i$.

We consider $M$ as an $S$-module via $\phi$. Since 
$M/(\X)M=M/(\x)M$ is a finite length $A$-module and so a 
finite length $S$-module, since $S/\q S \cong A/\m$. 
So $M$ is a finite $S$-module. Also note that
\[
\q = \sqrt{\ann_S(M/\X M)} = \sqrt{\ann_S(M) + (\X)}.
\]
So there exists $\Delta \in \ann_S(M) \setminus (\X)$. 
Observe that $\Delta,X_1,\ldots,X_n$ is an s.o.p. of $S$. Since $S$ is regular local ring of
dimension $n+1$, we have that $\Delta, X_1,\ldots,X_n$ is an $S$-regular 
sequence. 
Set $B=S/\Delta$ and $y_i=\overline{X}_i$ for $i=1,\ldots,n.$
Note that $B$ satisfies our requirements.
\end{proof}

\bigskip\noindent
{\bf Acknowledgment :} The author  thanks Prof. W. Bruns and Prof. J. Herzog  
for helpful discussions.

\providecommand{\bysame}{\leavevmode\hbox to3em{\hrulefill}\thinspace}
\providecommand{\MR}{\relax\ifhmode\unskip\space\fi MR }
\providecommand{\MRhref}[2]{%
  \href{http://www.ams.org/mathscinet-getitem?mr=#1}{#2}
}
\providecommand{\href}[2]{#2}

\end{document}